\documentclass[12pt]{scrartcl}

\usepackage{setspace}

\usepackage{color}
\usepackage{array}
\usepackage{supertabular}
\usepackage{hhline}
\usepackage{multirow}
\usepackage{multicol}
\usepackage{textgreek}
\usepackage{upgreek}

\usepackage{wrapfig}

\usepackage{mathbbol}

\usepackage[listings,theorems]{tcolorbox}

\usepackage{pifont}

\usepackage[thicklines]{cancel}

\usepackage[LGR,T1]{fontenc}

\usepackage{tikz}

\usetikzlibrary{hobby}

\usepackage{pgfplots}

\pgfplotsset{compat=1.10}
\usepgfplotslibrary{fillbetween}
\usetikzlibrary{patterns}

 \usetikzlibrary{patterns}

\usepackage[latin1]{inputenc}

\usepackage{graphicx}
\usepackage{fancybox}

\usepackage[latin1]{inputenc}

\newcommand\blfootnote[1]{%
  \begingroup
  \renewcommand\thefootnote{}\footnote{#1}%
  \addtocounter{footnote}{-1}%
  \endgroup
}

\usepackage{graphicx,fancybox,manfnt,wrapfig,multicol,mathrsfs}

\usepackage{amsmath,amssymb,color,makeidx,rotating,pifont}

\usepackage{mathtools}

\usepackage{stmaryrd}

\mathtoolsset{showonlyrefs=true}  

\begin{document}

\LARGE\sffamily 
\begin{center}\bfseries
On fractional differential equations,\\
dimensional analysis, and\\
the double gamma function{\blfootnote{Preprint version -- submitted for publication in Nonlinear Dynamics}}
\end{center}

\normalsize

\begin{center}
{\bfseries J. Vaz\footnote{\texttt{vaz@unicamp.br}} and E. Capelas de Oliveira\footnote{\texttt{capelas@unicamp.br}}}
\end{center}

\normalsize\rmfamily\sffamily

\begin{center}
Departamento de Matem\'atica Aplicada\\
Universidade Estadual de Campinas\\
Campinas, SP, Brazil
\end{center}

\bigskip

\small\rmfamily

\begin{center}
\begin{minipage}[c]{10cm}
\centerline{\bfseries\sffamily Abstract}
In this paper we discuss some issues that arise in the process of writing a 
fractional differential equation (FDE) 
by replacing an integer order derivative by 
a fractional order derivative in a given 
differential equation. To address these issues, we propose a dimensional regularization 
of the Caputo fractional derivative, ensuring consistency in physical dimensions. 
Then we solve some FDEs using this proposed dimensional regularization. 
We show that the solutions of these FDEs are most conveniently written using 
the double gamma function. We also compare these solutions with 
those from equations involving the standard Caputo fractional derivative. 
\end{minipage}
\end{center}

\medskip

\centerline{\textit{Dedicated to the memory of Jos\'e Ant\'onio Tenreiro Machado}}

\normalsize\rmfamily

\medskip

\DeclareGraphicsExtensions{.gif,.pdf}

\normalsize

\section{Introduction}

The motivation for this work comes from the problem of addressing the physical dimensions in a FDE. 
Indeed, from a purely mathematical 
point of view, considering a differential equation (DE) involving variables with 
physical dimensions is a completely unnecessary complication. However, 
when we want to use these DEs to model real problems, 
the variables involved in it must have physical dimensions. 
For DE involving (ordinary or partial) derivatives of integer order, 
dealing with problems involving variables with physical dimensions  
does not represent a problem from the mathematical point of view, after all
it is always possible to transform a DE involving variables with physical dimensions 
into a dimensionless form. This transformation relies on the process of 
\textit{nondimensionalization}, where appropriate dimensionless variables 
and parameters are introduced using the characteristic scales of the system.

However, we have some subtleties when we consider a DE involving non-integer order derivatives \cite{dim1,dim2}. 
In Section 2 we will discuss some examples where we will see that the use of a 
fractional derivative in a DE requires caution and can 
lead to problems. Indeed, let us consider the Caputo fractional differential
operator, which is widely used in FDE. If $x$ has physical dimension $\mathrm{L}$ 
and $t$ has physical dimension $\mathrm{T}$, then the physical dimensional of the
Caputo fractional derivative ${\sideset{_{}^{}}{_t^{\alpha}}{\operatorname{\mathcal{D}}}}x$ 
is $\mathrm{L} \mathrm{T}^{-\alpha}$. Thus, if in a DE involving variables with physical dimensions 
and a derivative, for example, $dx/dt$, we assume its fractional version as being 
the equation obtained by simply replacing $dx/dt$ with ${\sideset{_{}^{}}{_t^{\alpha}}{\operatorname{\mathcal{D}}}}x$ , the result will be an unbalanced equation from 
the point of view of physical dimensions. Thus, the question we must ask is: 
how to write a fractional version of a DE with integer-order derivatives with 
the correct dimensional balance? 

In addition to discussing the problems 
involving dimensional balancing in a FDE, we will also discuss how to make a 
small modification in the use of the Caputo fractional derivative, introducing 
a multiplicative term involving the independent variable  in order to maintain 
the same original physical dimension of the problem with integer order derivative. 
This modification does not change the main characteristic of the Caputo fractional 
derivative, which is its non-locality, and only corrects the problem of physical 
dimension. We will see that for DE of order two (or higher) 
there are two natural ways to make this modification. 

In Section 3 we solve two FDE models based on these modifications (which we 
will call Rule I and Rule II). An interesting thing about FDEs is 
that their study often leads to 
the introduction of new mathematical functions, or at least functions that are 
not very popular in the scientific literature. We will see that the solutions 
of these FDEs 
are given in terms of products of gamma functions, and to write these products 
in a compact form, we will show that it is convenient to use Barnes $G$-function $G(z)$ 
or its generalization, the double gamma function $G(z;\tau)$. In fact, some authors
also refers to $G(z)$ as double gamma function, which is not a mistake since
$G(z) = G(z;1)$. Therefore, unless we are not specifically restricting ourselves 
to the case $\tau = 1$, we will refer to $G(z)$ and $G(z;\tau)$ interchangeably as
double gamma functions. 
Despite generalizing a very important mathematical function, that is the gamma
function, the interest in these functions was not great, even though they appear as 
an exercise in Whittaker and Watson's famous book \cite{WW}. However, in recent decades, 
some applications have emerged in the literature. We can mention, for example,
their appearance in the study 
of the determinants of the Laplacians on the $n$-dimensional unit sphere 
$S^n$ \cite{Choi,Kumagai} and in a proof of the Kronecker limit formula \cite{Shintani}. 
A great advantage of using 
the double gamma function  instead of the
product of gamma functions is that these functions have integral representations that 
can be used to compute their numerical values. 
Since the literature on these functions is not large, we believe 
it would be useful for the reader to include an appendix on them 
(Appendix~A) as well as an appendix on Kilbas-Saigo functions (Appendix~B), 
which naturally arise in the solutions of one of these FDE models.

In section 4 we give a comparison between the functions that were found in the previous section with other functions that arise as solutions of other FDEs, in particular with the Mittag-Leffler function. 
In Section 5 we present our concluding remarks.

\section{FDE and physical dimensions}

There are many different definitions of fractional derivative in the 
literature \cite{Capelas,Capelas2}, but when we are interested in an initial 
value problem, certainly the best known and most relevant is the 
Caputo fractional derivative, that is, 
\begin{equation}
\label{caputo.def}
{\sideset{_{}^{}}{_t^{\alpha}}{\operatorname{\mathcal{D}}}} y = 
\frac{1}{\Gamma(N-\alpha)} 
\int_0^t \frac{ \phi^{(N)}(\tau)}{(t-\tau)^{\alpha+1-N}} \, d\tau ,
\end{equation}
where $y = \phi(t)$ satisfying appropriate conditions, 
$N -1 < \alpha \leq N$ ($N=1,2\ldots$), and the initialization was chosen at $a=0$. 
An important characteristic of the Caputo fractional derivative is the initial
conditions are the same in an initial value problem with Caputo fractional
derivative and ordinary integer-order derivative. 
Therefore, we will use the Caputo fractional derivative in what follows. 
Furthermore, to discuss the problem we have in mind, we believe the best approach 
is through examples, which is what we will do next, where we will  
denote variables with physical dimensions 
by $x, y, t, \ldots$, with $[x], [y], [t],\ldots$ being the respective physical dimensions, 
and by $\bar{x},\bar{y}, \bar{t},\ldots$ the respective dimensionless variables. 
We will also define characteristic quantities as, for example, $t_{\text{c}}$
such that $t = t_\text{c} \bar{t}$, and so on. 

\paragraph{Example 1.} Let us start our discussion with the decay or relaxation equation 
\begin{equation}
\label{eq.relax.2}
\frac{dx}{dt} + \kappa x = 0 .
\end{equation}
where $\kappa$ is a constant with physical dimension $[\kappa] = \mathrm{T}^{-1}$, where
$\mathrm{T} = [t]$. Its version in dimensionless form is 
\begin{equation}
\label{eq.relax.1}
\frac{d\bar{x}}{d\bar{t}} + \bar{\kappa} \bar{x} = 0 , 
\end{equation}
where $\bar{\kappa}$ is a positive constant. 
Eq.\eqref{eq.relax.1} and eq.\eqref{eq.relax.2} 
are related by the change of variables $t = t_{\text{c}} \bar{t}$ with $t_{\text{c}}= \bar{\kappa}/\kappa$ 
and $x = x_{\text{c}} \bar{x}$ with $x_{\text{c}}$ arbitrary. 

We want to consider a version of eq.\eqref{eq.relax.1} and 
eq.\eqref{eq.relax.2} in terms of the Caputo derivative of 
order $\alpha$ ($0 <\alpha \leq 1$). Since the Caputo derivative 
${\sideset{_{}^{}}{_{\bar{t}}^{\alpha}}{\operatorname{\mathcal{D}}}}\bar{x}$ 
is dimensionless, the fractional version of eq.\eqref{eq.relax.1} 
does not present any difficulty from the point of view 
of its dimensionality, and it is  
\begin{equation}
\label{eq.relax.fd1}
{\sideset{_{}^{}}{_{\bar{t}}^{\alpha}}{\operatorname{\mathcal{D}}}}\bar{x} + \bar{\kappa} \bar{x} = 0 . 
\end{equation}
However, the fractional version of eq.\eqref{eq.relax.2} requires a 
small adjustment since $[{\sideset{_{}^{}}{_{{t}}^{\alpha}}{\operatorname{\mathcal{D}}}} {x}] = \mathrm{T}^{-\alpha} \mathrm{L}$, 
where $\mathrm{T} = [t]$ and $\mathrm{L} = [x]$. The fractional version of eq.\eqref{eq.relax.2} 
can be obtained from eq.\eqref{eq.relax.fd1} using the relation between the variables, which gives  
\begin{equation}
\label{eq.relax.fd3.0}
{\sideset{_{}^{}}{_{t}^{\alpha}}{\operatorname{\mathcal{D}}}}x + 
\bar{\kappa}^{1-\alpha} \kappa^\alpha x = 0  
\end{equation}
or 
\begin{equation}
\label{eq.relax.fd3.0.1}
t_\text{c}^{\alpha-1} {\sideset{_{}^{}}{_{t}^{\alpha}}{\operatorname{\mathcal{D}}}}x + 
 \kappa  x = 0  .
\end{equation}
Note that the quantity $\bar{\kappa}^{1-\alpha}$ is a number, and we can always choose
$\bar{\kappa} = 1$ in eq.\eqref{eq.relax.fd3.0} or $t_\text{c} = 1/\kappa$ in eq.\eqref{eq.relax.fd3.0.1} 
since this is equivalent to use 
a change of scale $\bar{t} \to \lambda \bar{t}$ with $\lambda = 1/\bar{\kappa}$. 
Thus, without loss of generality, we can assume that the fractional version of eq.\eqref{eq.relax.2} is 
\begin{equation}
\label{eq.relax.fd3}
{\sideset{_{}^{}}{_{t}^{\alpha}}{\operatorname{\mathcal{D}}}}x +  \kappa^\alpha x = 0  .
\end{equation}

On the other hand, when we think about writing a new DE with the fractional derivative 
of order $\alpha$ generalizing a given DE with derivative of order 1, 
what we should expect is this new DE to be coherent from a dimensional point of view and 
that, in the limit $\alpha \to 1$, to reproduce the original equation. 
Within this perspective, we have another possibility, suggested by 
the form of eq.\eqref{eq.relax.fd3.0.1}. Indeed, since $[t_\text{c}^{\alpha-1}] = 
\mathrm{T}^{\alpha-1}$, this term re-establishes the dimensional balancing of the
equation. But another possibility is to use $t$ instead of $t_\text{c}$, that is,  
\begin{equation}
\label{eq.relax.fd4}
\bar{t}^{\alpha-1} {\sideset{_{}^{}}{_{\bar{t}}^{\alpha}}{\operatorname{\mathcal{D}}}}\bar{x} + \bar{\kappa} \bar{x} = 0   
\end{equation}
or 
\begin{equation}
\label{eq.relax.fd5}
t^{\alpha-1} {\sideset{_{}^{}}{_{t}^{\alpha}}{\operatorname{\mathcal{D}}}}x + \kappa x = 0  .
\end{equation}
Eq.\eqref{eq.relax.fd4} and eq.\eqref{eq.relax.fd5} are also related 
by $t = t_{\text{c}}  \bar{t} $ and $x = x_{\text{c}} \bar{x}$, like eq.\eqref{eq.relax.2} and 
eq.\eqref{eq.relax.1}. But unlike eq.\eqref{eq.relax.fd3} where we have
assumed $\bar{\kappa} = 1$, in eq.\eqref{eq.relax.fd5} we did not have 
to make any choice concerning $t_\text{c}$ or $\bar{\kappa}$. 

The fact is that both  eq.\eqref{eq.relax.fd3} 
and eq.\eqref{eq.relax.fd5} are two dimensionally 
coherent FDEs and are fractional generalizations of eq.\eqref{eq.relax.1}. 
At this point there is no a priori reason 
to give any preference to one of the two sets. Certainly, if we are thinking 
of using these equations to model anomalous relaxation processes, for example, 
then there is a criterion for giving a preference to one or the other set, 
which is the adequacy to the experimental results.  

\medskip

\paragraph{Example 2.} Let $x$ and $t$ have physical dimensions $\mathrm{L}$ and $\mathrm{T}$, 
respectively, and consider the DE 
\begin{equation}
\label{eq.ex.2.1}
\frac{dx}{dt} = a + b x + c x^2 ,
\end{equation}
where we assume that $a, b, c > 0$. The physical dimensions of $a$, $b$ and $c$ are
\begin{equation}
[a] = \mathrm{L}\mathrm{T}^{-1}, \quad [b] = \mathrm{T}^{-1}, \qquad [c] = \mathrm{L}^{-1} 
\mathrm{T}^{-1} . 
\end{equation}
Moreover, let us also consider the dimensionless DE 
\begin{equation}
\label{eq.ex.2.2}
\frac{d\bar{x}}{d\bar{t}} = \bar{a} + \bar{b} \bar{x} + \bar{c} \bar{x}^2 .
\end{equation}
If we write $t = t_{\text{c}} \bar{t}$ and $x = x_{\text{c}} \bar{x} $, with 
$[t_{\text{c}}] = \mathrm{T}$ and $[x_{\text{c}}] = \mathrm{L}$, we have
\begin{equation}
\label{eq.ex.2.3}
\bar{a} = \frac{t_{\text{c}} a}{x_{\text{c}}} , \quad 
\bar{b} = t_{\text{c}} b , \quad \bar{c} = t_{\text{c}} x_{\text{c}} c . 
\end{equation}
From these relations it follows that we can define a dimensionless parameter $\varkappa$ as 
\begin{equation}
\label{eq.ex.2.4}
\varkappa = \frac{ac}{b^2} = \frac{\bar{a}\bar{c}}{{\bar{b}}^2} . 
\end{equation}
As a consequence, unlike in the previous example where we can always choose the parameter 
$\bar{\kappa} = 1$ in eq.\eqref{eq.relax.1}, in eq.\eqref{eq.ex.2.2} we can only 
choose arbitrarily the values of any two parameters from $\bar{a}$, $\bar{b}$ and $\bar{c}$, 
the other one being determined by the relation in eq.\eqref{eq.ex.2.4} -- for example, 
$(\bar{a},\bar{b},\bar{c}) = \{(1,1,\varkappa),(\varkappa,1,1),(1,{\varkappa}^{-1/2},1)\}$ -- 
and this is due to the fact that we have three parameters but only two 
variables.

Again, since the Caputo derivative 
${\sideset{_{}^{}}{_{\bar{t}}^{\alpha}}{\operatorname{\mathcal{D}}}}\bar{x}$ 
is dimensionless, the fractional version of eq.\eqref{eq.ex.2.1} is 
\begin{equation}
\label{eq.ex.2.5}
{\sideset{_{}^{}}{_{\bar{t}}^{\alpha}}{\operatorname{\mathcal{D}}}}\bar{x} = 
\bar{a} + \bar{b} \bar{x} + \bar{c} \bar{x}^2 .
\end{equation}

For the dimensional version, using $t = t_{\text{c}} \bar{t}$, $x = x_{\text{c}} \bar{x} $ and eq.\eqref{eq.ex.2.3}, 
we obtain 
\begin{equation}
\label{eq.ex.2.6}
{\sideset{_{}^{}}{_{{t}}^{\alpha}}{\operatorname{\mathcal{D}}}}{x} = 
t_{\text{c}}^{1-\alpha}(a + b x + c x^2 ) .
\end{equation}
The situation is now a little more complicate because the dimensional factor $t_{\text{c}}^{1-\alpha}$ 
cannot be absorbed into powers of the terms of the equation, as in the previous example --  
indeed there is no $\beta$ such that $a^\beta$ has dimension $\mathrm{L}\mathrm{T}^{-1}$.  
On the other hand, recalling from eq.\eqref{eq.ex.2.3} that $t_{\text{c}} = \bar{b}/b$, we can also 
write 
\begin{equation}
\label{eq.ex.2.6.5}
{\sideset{_{}^{}}{_{{t}}^{\alpha}}{\operatorname{\mathcal{D}}}}{x} = 
\bar{b}^{1-\alpha}(a b^ {\alpha-1} + b^ \alpha x + c b^ {\alpha-1} x^2 ) . 
\end{equation}
The factor $\bar{b}^{1-\alpha}$ is dimensionless. So we have on the
right hand side a term 
$b^\alpha$ like $\kappa^ \alpha$ in the previous example, but the other 
ones involve also the parameter $b$, that is, $a b^{\alpha-1}$ and 
$c b^{\alpha-1}$. Note that this brings us to a problem concerning 
the interpretation of the equation. If the terms $a$, $b$ and $c$ 
have different physical interpretations, what is the interpretation 
of the terms $a b^{\alpha-1}$ and $c b^{\alpha-1}$ involving 
combinations of terms with different individual interpretations? 

Like in the previous example, we can always choose an appropriate scale such that $\bar{b}=1$ (for 
example, related to the choices $\{(1,1,\varkappa),(\varkappa,1,1)\})$ for $(\bar{a},\bar{b},\bar{c})$, so, 
without loss of generality, we can assume the fractional version of 
eq.\eqref{eq.ex.2.1} as 
\begin{equation}
\label{eq.ex.2.6.6}
{\sideset{_{}^{}}{_{{t}}^{\alpha}}{\operatorname{\mathcal{D}}}}{x} = 
a b^ {\alpha-1} + b^ \alpha x + c b^ {\alpha-1} x^2  . 
\end{equation}

\medskip

There is still another possibility, as seen in the previous example. Let us
rewrite eq.\eqref{eq.ex.2.6} as 
\begin{equation}
\label{eq.ex.2.6.alt}
t_{\text{c}}^{\alpha-1} {\sideset{_{}^{}}{_{{t}}^{\alpha}}{\operatorname{\mathcal{D}}}}{x} = 
(a + b x + c x^2 ) .
\end{equation}
Instead of $t_\text{c}$, we can consider the variable $t$, and consequently 
the equations
\begin{equation}
\label{eq.ex.2.7}
\bar{t}^{\alpha-1}{\sideset{_{}^{}}{_{\bar{t}}^{\alpha}}{\operatorname{\mathcal{D}}}}\bar{x} = 
\bar{a} + \bar{b} \bar{x} + \bar{c} \bar{x}^2  
\end{equation}
and 
\begin{equation}
\label{eq.ex.2.8}
t^{\alpha-1} {\sideset{_{}^{}}{_{{t}}^{\alpha}}{\operatorname{\mathcal{D}}}}{x} =  a + b x + c x^2 . 
\end{equation}
Note that eq.\eqref{eq.ex.2.8} follows from eq.\eqref{eq.ex.2.7} using $t = t_{\text{c}} \bar{t}$, $x = x_{\text{c}} \bar{x} $ 
and eq.\eqref{eq.ex.2.3}. Moreover, the terms on the right hand side have the same
interpretation both in the first-order and the fractional order derivatives.

\medskip

\paragraph{Example 3.} An important example is the damped harmonic oscillator equation, 
\begin{equation}
\label{eq.ex.3.1}
\frac{d^2 x}{dt^2} + 2\gamma \frac{dx}{dt} + \omega^2 x = 0 ,
\end{equation}
where $m = 1$ for simplicity 
and the physical dimension of the damping constant $\gamma$ and the frequency $\omega$ are
\begin{equation}
[\gamma] = [\omega] = \mathrm{T}^{-1} . 
\end{equation}
Its dimensionless form is 
\begin{equation}
\label{eq.ex.3.2}
\frac{d^2 \bar{x}}{d\bar{t}^2} + 2\bar{\gamma} \frac{d\bar{x}}{d\bar{t}} + \bar{\omega}^2 \bar{x} = 0 , 
\end{equation}
where $t = t_{\text{c}} \bar{t}$, $x = x_{\text{c}} \bar{x}$ and
\begin{equation}
\label{eq.ex.3.3.0}
\bar{\gamma} = t_{\text{c}} \gamma, \qquad \bar{\omega} = t_{\text{c}} \omega . 
\end{equation} 

The fractional version of eq.\eqref{eq.ex.3.2} is 
\begin{equation}
\label{eq.ex.3.3}
{\sideset{_{}^{}}{_{\bar{t}}^{\beta}}{\operatorname{\mathcal{D}}}}\bar{x} + 
2\bar{\gamma} {\sideset{_{}^{}}{_{\bar{t}}^{\alpha}}{\operatorname{\mathcal{D}}}}\bar{x} 
+ \bar{\omega}^2 \bar{x} = 0 , 
\end{equation}
where $1 < \beta \leq 2$ and $0 < \alpha \leq 1$. Returning to the dimensional variables, we 
have 
\begin{equation}
\label{eq.ex.3.4}
{\sideset{_{}^{}}{_{{t}}^{\beta}}{\operatorname{\mathcal{D}}}} {x} + 
2 {\gamma} t_{\text{c}}^{1+\alpha-\beta} {\sideset{_{}^{}}{_{{t}}^{\alpha}}{\operatorname{\mathcal{D}}}} {x} 
+  {\omega}^2  t_{\text{c}}^{2-\beta} {x} = 0 , 
\end{equation}
where $t_{\text{c}}$ is an arbitrary dimensional parameter. Using eq.\eqref{eq.ex.3.3.0}, 
we can write this equation as 
\begin{equation}
\label{eq.ex.3.4.1}
{\sideset{_{}^{}}{_{{t}}^{\beta}}{\operatorname{\mathcal{D}}}} {x} + 
2 {\gamma}^{\beta-\alpha} \bar{\gamma}^{1-\alpha-\beta} {\sideset{_{}^{}}{_{{t}}^{\alpha}}{\operatorname{\mathcal{D}}}} {x} 
+  {\omega}^2  \gamma^{\beta-2} \bar{\gamma}^{2-\beta} {x} = 0 , 
\end{equation}
or 
\begin{equation}
\label{eq.ex.3.4.2}
{\sideset{_{}^{}}{_{{t}}^{\beta}}{\operatorname{\mathcal{D}}}} {x} + 
2 {\gamma} \omega^{\beta-\alpha-1} 
\bar{\omega}^{1+\alpha-\beta} {\sideset{_{}^{}}{_{{t}}^{\alpha}}{\operatorname{\mathcal{D}}}} {x} 
+  {\omega}^\beta  \bar{\omega}^{2-\beta} {x} = 0 , 
\end{equation}
where $\bar{\gamma}$ and $\bar{\omega}$ are dimensionless. 
Note that we also have the problem of interpreting the terms 
of the above equations. In eq.\eqref{eq.ex.3.4.1} we have
a frequency-like term ${\omega}^2  \gamma^{\beta-2}$ 
involving the damping factor $\gamma$ when $\beta \neq 2$ and 
in eq.\eqref{eq.ex.3.4.2} we have a damping-like term 
$2 {\gamma} \omega^{\beta-\alpha-1} $ involving the
frequency $\omega$ when $\beta -\alpha \neq 1$.

We can choose a specific value for $\bar{\gamma}$ or 
$\bar{\omega}$ to simplify the equations. The most interesting choice is 
$\bar{\omega} = 1$. In this case $t_{\text{c}} = \omega^{-1}$ and eq.\eqref{eq.ex.3.4} becomes 
\begin{equation}
\label{eq.ex.3.5}
{\sideset{_{}^{}}{_{{t}}^{\beta}}{\operatorname{\mathcal{D}}}} {x} + 
2 {\gamma} \omega^{\beta-\alpha-1} {\sideset{_{}^{}}{_{{t}}^{\alpha}}{\operatorname{\mathcal{D}}}} {x} 
+  {\omega}^\beta    {x} = 0 . 
\end{equation}
Moreover, when $\beta=\alpha-1$, we have 
\begin{equation}
\label{eq.ex.3.6}
{\sideset{_{}^{}}{_{{t}}^{\beta}}{\operatorname{\mathcal{D}}}} {x} + 
2 {\gamma}   {\sideset{_{}^{}}{_{{t}}^{\beta-1}}{\operatorname{\mathcal{D}}}} {x} 
+  {\omega}^\beta    {x} = 0 . 
\end{equation}
When $\gamma = 0$ we have 
 \begin{equation}
{\sideset{_{}^{}}{_{{t}}^{\beta}}{\operatorname{\mathcal{D}}}} {x} + 
  {\omega}^\beta    {x} = 0 ,
\end{equation}
which resembles eq.\eqref{eq.relax.fd3},  
except that here $1 < \beta \leq 2$ while $0 < \alpha \leq 1$ in eq.\eqref{eq.relax.fd3}.

For the choice $\bar{\gamma} = 1$, that is, $t_\text{c} = \gamma^{-1}$,  we obtain, instead of eq.\eqref{eq.ex.3.5}, 
\begin{equation}
{\sideset{_{}^{}}{_{{t}}^{\beta}}{\operatorname{\mathcal{D}}}} {x} + 
2 {\gamma}^{\beta-\alpha} {\sideset{_{}^{}}{_{{t}}^{\alpha}}{\operatorname{\mathcal{D}}}} {x} 
+  {\omega}^2  \gamma^{\beta-2} {x} = 0 .
\end{equation}

In a similar manner to the previous examples, there is an alternative 
version for the above equations, which are 
\begin{equation}
\label{eq.ex.3.7}
t^{\beta-2} {\sideset{_{}^{}}{_{\bar{t}}^{\beta}}{\operatorname{\mathcal{D}}}}\bar{x} + 
2\bar{\gamma} t^{\alpha-1} {\sideset{_{}^{}}{_{\bar{t}}^{\alpha}}{\operatorname{\mathcal{D}}}}\bar{x} 
+ \bar{\omega}^2 \bar{x} = 0 , 
\end{equation}
and 
\begin{equation}
\label{eq.ex.3.8}
t^{\beta-2} {\sideset{_{}^{}}{_{{t}}^{\beta}}{\operatorname{\mathcal{D}}}}{x} + 
2{\gamma} t^{\alpha-1} {\sideset{_{}^{}}{_{{t}}^{\alpha}}{\operatorname{\mathcal{D}}}}{x} 
+ {\omega}^2 {x} = 0 , 
\end{equation}
with $1 < \beta \leq 2$ and $0 <\alpha \leq 1$. 
Eq.\eqref{eq.ex.3.7} and eq.\eqref{eq.ex.3.8} are related by $t = t_{\text{c}}\bar{t}$ and 
$x = x_{\text{c}} \bar{x}$ and eq.\eqref{eq.ex.3.3.0}. 
Moreover, there is an alternative to eq.\eqref{eq.ex.3.7} and eq.\eqref{eq.ex.3.8} which is also 
consistent from the dimensional point of view, that is 
\begin{equation}
\label{eq.ex.3.9}
\left(t^{\alpha^\prime-1} {\sideset{_{}^{}}{_{\bar{t}}^{\alpha^\prime}}{\operatorname{\mathcal{D}}}}\right)^2 \bar{x} + 
2\bar{\gamma} t^{\alpha-1} {\sideset{_{}^{}}{_{\bar{t}}^{\alpha}}{\operatorname{\mathcal{D}}}}\bar{x} 
+ \bar{\omega}^2 \bar{x} = 0 , 
\end{equation}
and 
\begin{equation}
\label{eq.ex.3.10}
\left(t^{\alpha^\prime-1} {\sideset{_{}^{}}{_{{t}}^{\alpha^\prime}}{\operatorname{\mathcal{D}}}}\right)^2{x} + 
2{\gamma} t^{\alpha-1} {\sideset{_{}^{}}{_{{t}}^{\alpha}}{\operatorname{\mathcal{D}}}}{x} 
+ {\omega}^2 {x} = 0 , 
\end{equation}
where $0 < \alpha, \alpha^\prime \leq 1$. 

\bigskip

\paragraph{Remark.} The above examples have shown that simply replacing 
an integer-order derivative by a Caputo fractional order derivative 
in a DE is a procedure that leads to a well-defined FDE 
\textit{in the case of dimensionless variables}. However, when the 
DE involves variables with physical dimensions, this procedure 
has encountered difficulties, particularly in the issue 
of the \textit{dimensional consistency} of the equation and 
the \textit{interpretation} of its terms.

On the other hand, we saw in the examples 
that this problem of dimensional balancing of the equation with Caputo 
fractional derivative can be solved with the help of the independent variable  
involved in the problem. We discussed two ways to do this: one is doing 
the replacement according to the rule 
\begin{equation}
\label{subs.1}
\text{(I)} \qquad \frac{d^N\;}{d t^N} \to 
 t^{\alpha - N} {\sideset{_{}^{}}{_{{t}}^{\alpha}}{\operatorname{\mathcal{D}}}} 
\end{equation}
with $N-1 < \alpha \leq N$ ($N=1,2,\ldots$), 
and the other is doing the replacement according to the rule 
\begin{equation}
\label{subs.2}
\text{(II)} \qquad \frac{d^N\;}{d t^N} \to  \left(t^{\alpha^\prime - 1} {\sideset{_{}^{}}{_{{t}}^{\alpha^\prime}}{\operatorname{\mathcal{D}}}}\right)^N 
\end{equation}
with $0 < \alpha^\prime \leq 1$ ($N = 1,2,\ldots$). Note that for $[t] = \mathrm{T}$ we have 
\begin{equation}
\label{aux.dim}
 \left[\displaystyle{\frac{d^N\;}{dt^N}}\right] = [ t^{\alpha - N} {\sideset{_{}^{}}{_{{t}}^{\alpha}}{\operatorname{\mathcal{D}}}} ] = 
\left[\left(t^{\alpha^\prime - 1} {\sideset{_{}^{}}{_{{t}}^{\alpha^\prime}}{\operatorname{\mathcal{D}}}}\right)^N \right] = 
\mathrm{T}^{-N} . 
\end{equation}
In this way, we obtain a dimensionally consistent FDE without 
the need for major modifications to the equation and the 
interpretation of its terms. 
The fact that the physical dimension of the quantities in eq.\eqref{aux.dim} 
does not depend on $\alpha$ suggests the interpretation of $\alpha$, or better $N -\alpha$, 
as a measure of the \textit{degree of non-locality} involved in the model.

\section{FDE with dimensional regularized Caputo derivatives and double gamma functions}

Given that the FDE obtained from a DE with 
integer-order derivative through the rules in eq.\eqref{subs.1} or 
in eq.\eqref{subs.2} presents the same dimensional consistency as the original equation, 
the next step is to compare these FDEs, and one way to do this is by studying their solutions. 
In what follows we will study the fractional versions of the integer-order DE 
\begin{equation}
\label{DE.int.order.nu}
\frac{d^N f}{dt^N} + \omega^N f = 0 \qquad (N = 1,2,\ldots) 
\end{equation}
using the rules in eq.\eqref{subs.1} and eq.\eqref{subs.2}.  
Clearly $[\omega] = \mathrm{T}^{-1}$. 
This simple equation is enough to illustrate the difference between the 
two approaches and to illustrate the usefulness of the double gamma function 
in fractional calculus.

In what follows we will use 
\begin{equation}
\label{reg.caputo.n}
 t^{\alpha-N} {\sideset{_{}^{}}{_t^{\alpha}}{\operatorname{\mathcal{D}}}} t^m = 
\begin{cases}
{\displaystyle \frac{\Gamma(m+1)}{\Gamma(m-\alpha+1)} t^{m-N}} , \; & m > 0, \; m \notin \{0,1,2,\ldots,N-1\}, \\[1ex]
0 , & m \in \{0,1,\ldots,N-1\} ,
\end{cases}
\end{equation} 
for $N-1 < \alpha \leq N$ ($N = 1,2,\ldots$), which follows directly from eq.\eqref{caputo.def}. Thus, for 
$f(t)$ given by 
\begin{equation}
\label{analytic.1}
f(t) = \sum_{k=0}^\infty f_k t^k ,
\end{equation}
we have
\begin{equation}
\label{analytic.2}
t^{\alpha-N} {\sideset{_{}^{}}{_t^{\alpha}}{\operatorname{\mathcal{D}}}}f =
 \sum_{k=0}^\infty f_{k+N} \frac{\Gamma(k+N+1)}{\Gamma(k+N+1-\alpha)} t^k. 
\end{equation}

\subsection{FDE using rule I}

The fractional version of eq.\eqref{DE.int.order.nu} using the rule in eq.\eqref{subs.1} is  
\begin{equation}
\label{fde.rule.1}
t^{\alpha - N} {\sideset{_{}^{}}{_{{t}}^{\alpha}}{\operatorname{\mathcal{D}}}}  f = -\omega^N f , 
\end{equation}
where  $N-1 < \alpha \leq N$ for $N = 1,2,\ldots$ or $N = \lceil \alpha 
\rceil$, where $\lceil \cdot \rceil$ is the ceiling function. We will look for analytic solutions 
of eq.\eqref{fde.rule.1}. Using eq.\eqref{analytic.1} and eq.\eqref{analytic.2} in eq.\eqref{fde.rule.1} 
we obtain that the coefficients must satisfy 
\begin{equation}
\label{fde.1.coeff}
f_{k+N} \frac{\Gamma(k+N+1)}{\Gamma(k+N+1-\alpha)} = - \omega^N f_k , \qquad 
k = 0,1,2,\ldots 
\end{equation}
From this recurrence relation we obtain  
\begin{equation}
 f_{mN + j}  =  \left(-\omega^N\right)^m \left(\prod_{r=1}^m \frac{\Gamma(1+rN+j-\alpha)}{\Gamma(1+rN+j)}\right) f_j ,  
\end{equation}
with $j = 0,1,\ldots,N-1$ and $m = 1,2,\ldots$.  
Defining $C_{j,m}^{\scriptscriptstyle N}(\alpha)$ as 
\begin{equation}
\label{def.C}
C_{j,m}^{\scriptscriptstyle N}(\alpha) = \begin{cases} 1 , \quad & m = 0 , \\[1ex]
\displaystyle \prod_{r=1}^m  \frac{\Gamma(1+rN+j-\alpha)}{\Gamma(1+rN+j)} , \quad & m = 1,2,\ldots 
\end{cases}
\end{equation}
and the constants $\kappa_j$ as 
\begin{equation}
f_j = (\omega)^j \kappa_j ,
\end{equation}
we can write $f(t)$ as 
\begin{equation}
f(t) = \sum_{j=0}^{N-1} \kappa_j \phi_{j,\alpha}^{\scriptscriptstyle N}(\omega t) ,
\end{equation}
where 
\begin{equation}
\label{def_phi_j}
\phi_{j,\alpha}^{\scriptscriptstyle N}(\omega t) = \sum_{m=0}^\infty (-1)^m C_{j,m}^{\scriptscriptstyle N}(\alpha) (\omega t)^{m N + j} , \quad j = 0,1,\ldots,N-1.  
\end{equation}

The coefficients $C_{j,m}^{\scriptscriptstyle N}(\alpha)$ can be conveniently written 
using $G$-Pochhammer symbol (see Appendix~\ref{appendix.A}) as 
\begin{equation}
\label{coeff.G} 
C_{j,m}^{\scriptscriptstyle N}(\alpha) = \frac{\llbracket 1+(1+j-\alpha)\tau;\tau\rrbracket_m}{\llbracket 1+(1+j)\tau;\tau\rrbracket_m} 
\end{equation}
with 
\begin{equation}
\label{def.tau}
\tau = \frac{1}{N} , 
\end{equation}
where we recall that $N = \lceil \alpha \rceil = 1,2,\ldots$. Note that
eq.\eqref{coeff.G} also holds for $n=0$. 

We can also write 
\begin{equation}
\phi_{j,\alpha}^{\scriptscriptstyle N}(\omega t) = (\omega t)^j \sum_{m=0}^\infty   \frac{\llbracket 1+(1+j-\alpha)\tau;\tau\rrbracket_m}{\llbracket 1+(1+j)\tau;\tau\rrbracket_m}  \left(-(\omega t)^N\right)^{m} , 
\end{equation}
and comparing this series with the Kilbas-Saigo function (see Appendix B) as in eq.\eqref{KS.G} we conclude that 
\begin{equation}
\phi_{j,\alpha}^{\scriptscriptstyle N}(\omega t) =  (\omega t)^j E_{\alpha,\frac{N}{\alpha},\frac{N-\alpha+j}{\alpha}}\left(-(\omega t)^N\right) , 
\quad j = 0,1,\ldots,N-1 . 
\end{equation}

The convergence of $\phi_{j,\alpha}^{\scriptscriptstyle N}$ in eq.\eqref{def_phi_j} is guaranteed by Gautschi's inequality \cite{Gautschi}, which can be written as 
\begin{equation}
\frac{1}{x^{\sigma}} \leq \frac{\Gamma(x)}{\Gamma(x+\sigma)} \leq \left(1+\frac{1}{x}\right) 
\frac{1}{(x+1)^\sigma} ,
\end{equation}
for $x > 0$ and $0 \leq \sigma \leq 1$. In fact, using $x = 1 + (m+1)N + j-\alpha$ and 
$\sigma = \alpha$ it follows from this inequality and eq.\eqref{def.C} that 
\begin{equation}
\lim_{m\to \infty}\left|\frac{C_{j,m+1}(\alpha)}{C_{j,m}(\alpha)}\right| = 
\lim_{m\to \infty} \frac{\Gamma(1+(m+1)N+j-\alpha)}{\Gamma(1+(m+1)N+j)} = 0 , 
\end{equation}
and so eq.\eqref{def_phi_j} converges for all values of $\omega t$ and $0 < \alpha \leq 1$. 
The convergence for $\alpha$ such that $N-1 < \alpha \leq N$ with $N = 2,3,\ldots$ 
is a trivial consequence. 
 
\subsubsection{Examples} 

\begin{figure}[hbt]
\begin{center}
\includegraphics[width=12cm]{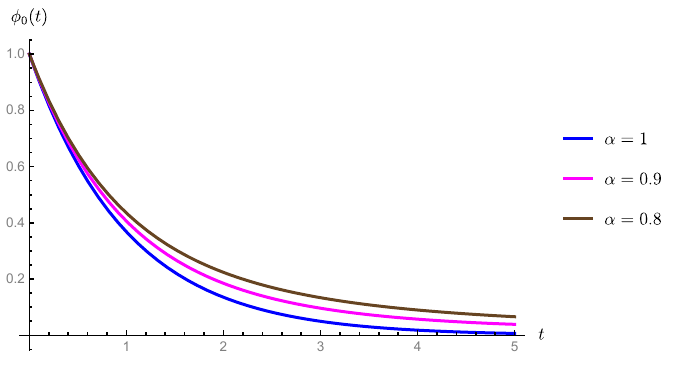}
\end{center}
\caption{Plots of $\phi_{0,\alpha}^{\scriptscriptstyle 1}(t)$ in eq.\eqref{phi.zero.alpha} for $\alpha = \{1,0.9,0.8\}$.
\label{fig.1}}
\end{figure}

\paragraph{(i) $\boldsymbol{N=1}$.} For $0 < \alpha \leq 1$ the solution of eq.\eqref{fde.rule.1} is 
\begin{equation}
\label{phi.zero.alpha}
\phi_{0,\alpha}^{\scriptscriptstyle 1}(\omega t) = \sum_{m=0}^\infty (-1)^m C_{0,m}^{\scriptscriptstyle 1}(\alpha) (\omega t)^{m} , 
\end{equation}
with 
 \begin{equation}
C_{0,m}^{\scriptscriptstyle 1}(\alpha) = \frac{\llbracket 2-\alpha;1\rrbracket_m}{\llbracket 2;1\rrbracket_m} = 
\frac{G(2-\alpha+m) G(2)}{G(2-\alpha) G(2+m)} , 
\end{equation}
where $G(\cdot)$ is the Barnes $G$-function. In Figure~\ref{fig.1} we have the
plots of $\phi_{0,\alpha}(t)$ for $\alpha = \{1,0.9,0.8\}$. 
Note that when $\alpha =1$ we have 
\begin{equation}
C_{0,m}^{\scriptscriptstyle 1}(1) =  
\frac{G(1+m) G(2)}{G(1) G(2+m)} = \frac{G(1+m) \Gamma(1) G(1)}{G(1)\Gamma(1+m)G(1+m)} = 
\frac{1}{\Gamma(1+m)} ,
\end{equation}
and so $\phi_{0,1}^{\scriptscriptstyle 1}(\omega t) = \operatorname{e}^{-\omega t}$. 

\begin{figure}[hbt]
\begin{center}
\includegraphics[width=12cm]{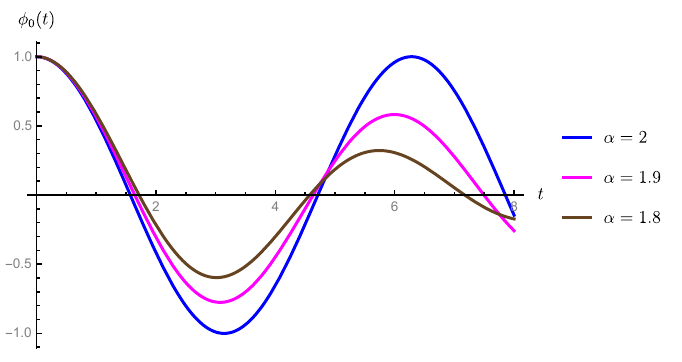}
\caption{Plots of $\phi_{0,\alpha}^{\scriptstyle 2}(t)$ in eq.\eqref{phi.zero.alpha.2} for $\alpha = \{2,1.9,1.8\}$.\label{fig.2}}
\includegraphics[width=12cm]{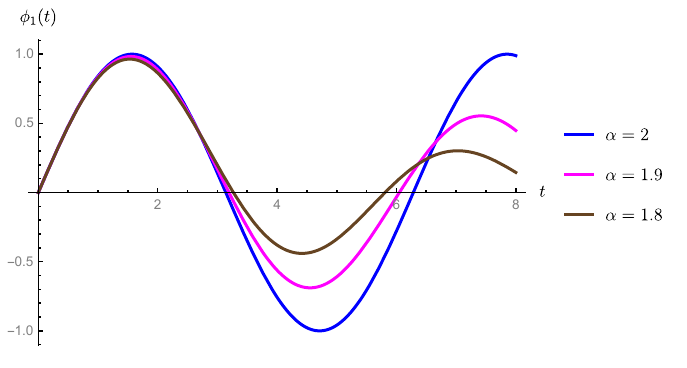}
\end{center}
\caption{Plots of $\phi_{1,\alpha}^{\scriptstyle 2}(t)$ in eq.\eqref{phi.one.alpha.2} for $\alpha = \{2,1.9,1.8\}$.\label{fig.3}}
 \end{figure}

\paragraph{(ii) $\boldsymbol{N=2}$.} For $1 < \alpha \leq 2$ the solutions of eq.\eqref{fde.rule.1}  
are 
\begin{equation}
\label{phi.zero.alpha.2}
\phi^{\scriptstyle 2}_{0,\alpha}(\omega t) = \sum_{m=0}^\infty (-1)^m C_{0,m}^{\scriptstyle 2}(\alpha) (\omega t)^{2m} ,
\end{equation}
with 
\begin{equation}
\label{C.0.m.aux}
C_{0,m}^{\scriptstyle 2}(\alpha) = \frac{\llbracket (3-\alpha)/2;1/2\rrbracket_m}{\llbracket 
3/2;1/2\rrbracket_m} = \frac{G[(3-\alpha)/2+m;1/2]G(3/2;1/2)}{G[(3-\alpha)/2;1/2]G(3/2+m;1/2)} , 
\end{equation}
and 
\begin{equation}
\label{phi.one.alpha.2}
\phi^{\scriptstyle 2}_{1,\alpha}(\omega t) = \sum_{m=0}^\infty (-1)^m C_{1,m}^{\scriptstyle 2}(\alpha) (\omega t)^{2m+1} ,
\end{equation}
with 
\begin{equation}
\label{C.0.m.aux.1}
C_{1,m}^{\scriptstyle 2}(\alpha) = \frac{\llbracket 2-\alpha/2;1/2\rrbracket_m}{\llbracket 
2;1/2\rrbracket_m} = \frac{G(2-\alpha/2+m;1/2)G(2;1/2)}{G(2-\alpha/2;1/2)G(2+m;1/2)} , 
\end{equation}
where $G(-;-)$ is the double gamma function. 
We have in Figure~\ref{fig.2} the
plots of $\phi_{0,\alpha}^{\scriptstyle 2}(t)$   
and in Figure~\ref{fig.3} the
plots of $\phi_{1,\alpha}^{\scriptstyle 2}(t)$ for $\alpha = \{2,1.9,1.8\}$. 
When $\alpha = 2$ we have, using eq.\eqref{ap.B.general.G.tau} in eq.\eqref{C.0.m.aux} 
and in  eq.\eqref{C.0.m.aux.1} that 
\begin{equation}
C_{0,m}^2(2) = \frac{G(1/2+m;1/2) G(3/2;1/2)}{G(1/2;1/2)G(3/2+m;1/2)} = \frac{1}{\Gamma(2m+1)}
\end{equation}
and 
 \begin{equation}
C_{1,m}^2(2) =  \frac{G(1+m;1/2) G(2;1/2)}{G(1;1/2)G(2+m;1/2)} = \frac{1}{\Gamma(2m+2)} ,
\end{equation}
and therefore 
\begin{equation}
\phi_{0,2}^2(\omega t) = \cos{\omega t} , \qquad \phi_{1,2}^2(\omega t) = \sin{\omega t} . 
\end{equation}

It is also interesting to note that using eq.\eqref{modular.2} with $p = \tau = 1$ and $q = 2$ 
we can rewrite $C_{j,m}^2(\alpha)$ for $j = 0,1$ as 
\begin{equation}
C_{j,m}^2(\alpha) = \frac{1}{2^{\alpha m}} \frac{\llbracket (3-\alpha+j)/2\rrbracket_m \llbracket (4-\alpha+j)/2\rrbracket_m}{\llbracket (3+j)/2\rrbracket_m \llbracket (4+j)/2\rrbracket_m} . 
\end{equation}

\subsection{FDE using rule II}

The fractional version of eq.\eqref{DE.int.order.nu} using the rule in eq.\eqref{subs.2} is  
\begin{equation}
\label{fde.rule.2}
\left(t^{\alpha^\prime - 1} {\sideset{_{}^{}}{_{{t}}^{\alpha^\prime}}{\operatorname{\mathcal{D}}}}\right)^N f = -\omega^N f , 
\end{equation}
where  $0 < \alpha^\prime \leq 1$ for $N = 1,2,\ldots$. In this case we have 
to use eq.\eqref{analytic.2} recursively, that is, 
\begin{equation}
\left(t^{\alpha^\prime - 1} {\sideset{_{}^{}}{_{{t}}^{\alpha^\prime}}{\operatorname{\mathcal{D}}}}\right)^N f = \sum_{k=0}^\infty f_{k+N} \left[\prod_{i=1}^N \frac{\Gamma(k+i+1)}{\Gamma(k+i+1-\alpha)}\right] t^k ,
\end{equation}
which in eq.\eqref{fde.rule.2} shows that we must have 
\begin{equation}
f_{k+N} = -\omega^N \prod_{i=1}^k \frac{\Gamma(k+i+1-\alpha^\prime)}{\Gamma(k+i+1)} f_k , \quad k=0,1,2,\ldots 
\end{equation}
From this recurrence relation we obtain 
\begin{equation}
f_{mN+j}  = \left(-\omega^N\right)^m \prod_{r=1}^{mN}\frac{\Gamma(r+j+1-\alpha^\prime)}{\Gamma(r+j+1)} f_j ,
\end{equation}
with $j = 0,1,\ldots,N-1$ and $ m=1,2,\ldots$.

Defining $K^{\scriptscriptstyle N}_{j,m}(\alpha^\prime)$ as 
\begin{equation}
K^{\scriptscriptstyle N}_{j,m}(\alpha^\prime) = \begin{cases} 1 , \quad & m = 0 , \\[1ex]
\displaystyle \prod_{r=1}^{mN}  \frac{\Gamma(r+j+1-\alpha^\prime)}{\Gamma(r+j+1)} , \quad & m = 1,2,\ldots 
\end{cases}
\end{equation}
and the constants $\kappa_j$ as 
\begin{equation}
f_j = (\omega)^j \kappa_j ,
\end{equation}
we can write $f(t)$ as 
\begin{equation}
f(t) = \sum_{j=0}^{N-1} \kappa_j \psi_{j,\alpha^\prime}^{\scriptscriptstyle N}(\omega t) ,
\end{equation}
where 
\begin{equation}
\psi_{j,\alpha^\prime}^{\scriptscriptstyle N}(\omega t) = \sum_{m=0}^\infty (-1)^m K^{\scriptscriptstyle N}_{j,m}(\alpha^\prime) (\omega t)^{m N + j} , \quad j = 0,1,\ldots,N-1.  
\end{equation}

Like in the previous case, we can also use the double gamma function to 
write a convenient expression for the coefficients $K^{\scriptscriptstyle N}_{j,m}(\alpha^\prime)$, 
but in this case the parameter $\tau$ is $\tau = 1$, that is, we have the
Barnes $G$-function $G(z)$. We have 
\begin{equation}
K^{\scriptscriptstyle N}_{j,m}(\alpha^\prime) = 
\frac{\llbracket j + 2-\alpha^\prime\rrbracket_{mN}}{\llbracket j+2\rrbracket_{mN}} ,
\end{equation}
and then 
\begin{equation}
\label{def_psi_j}
\psi_{j,\alpha^\prime}^{\scriptscriptstyle N}(\omega t) = \sum_{m=0}^\infty (-1)^m \frac{\llbracket j + 2-\alpha^\prime\rrbracket_{mN}}{\llbracket j+2\rrbracket_{mN}} (\omega t)^{m N + j} , \quad j = 0,1,\ldots,N-1.
\end{equation}

Like the case of eq.\eqref{def_phi_j}, the convergence of eq.\eqref{def_psi_j}
 for all values of $\omega t$ also follows using 
Gautschi's inequality. 

\begin{figure}[htb]
\begin{center}
\includegraphics[width=12cm]{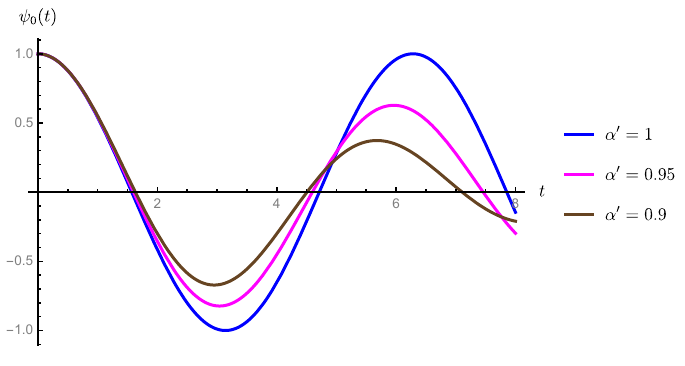}
\caption{Plots of $\psi_{0,\alpha^\prime}^{\scriptstyle 2}(t)$ in eq.\eqref{psi.0.alphaprime} for $\alpha^\prime = \{1,0.95,0.9\}$.\label{fig.5}}
\includegraphics[width=12cm]{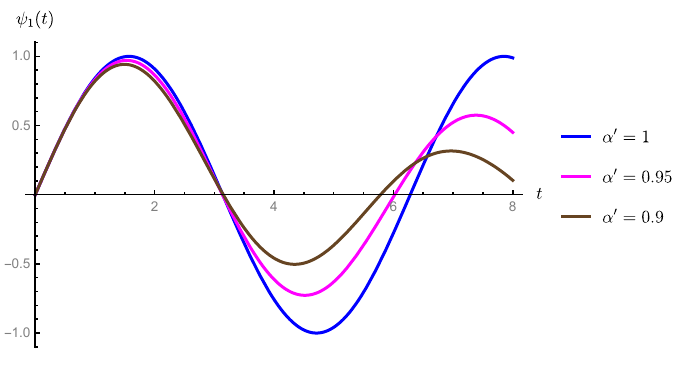}
\end{center}
\caption{Plots of $\psi_{1,\alpha^\prime}^{\scriptstyle 2}(t)$ in eq.\eqref{psi.1.alphaprime} for $\alpha^\prime = \{1,0.95,0.9\}$.\label{fig.6}}
 \end{figure}
 
\subsubsection{Examples}

\paragraph{(i) $\boldsymbol{N=2}$.} From eq.\eqref{def_psi_j} we have 
\begin{equation}
\label{psi.0.alphaprime}
\begin{aligned}
\psi_{0,\alpha^\prime}^2(\omega t) & = \sum_{m=0}^\infty \frac{\llbracket 2-\alpha^\prime\rrbracket_{2m}}{\llbracket 2 \rrbracket_{2m}} (-1)^m (\omega t)^{2m}   \\[1ex]
& = \frac{G(2)}{G(2-\alpha^\prime)} 
\sum_{m=0}^\infty \frac{G(2-\alpha^\prime+2m)}{G(2+2m)} (-1)^m (\omega t)^{2m}
\end{aligned}
\end{equation}
and 
\begin{equation}
\label{psi.1.alphaprime}
\begin{aligned}
\psi_{1,\alpha^\prime}^2(\omega t) & = \sum_{m=0}^\infty \frac{\llbracket 3-\alpha^\prime\rrbracket_{2m}}{\llbracket 3 \rrbracket_{2m}} (-1)^m (\omega t)^{2m+1} \\[1ex]
&  = \frac{G(3)}{G(3-\alpha^\prime)} 
\sum_{m=0}^\infty \frac{G(3-\alpha^\prime+2m)}{G(3+2m)} (-1)^m (\omega t)^{2m+1} .
\end{aligned}
\end{equation}
Using eq.\eqref{functional.barnes.G} it is not difficult to see that, for $\alpha^\prime \to 1$, 
we have, as expected, 
\begin{equation}
\psi_{0,1}^2(\omega t) = \cos{\omega t} , \qquad 
\psi_{1,1}^2(\omega t) = \sin{\omega t} .
\end{equation}
We have in Figure~\ref{fig.5} the
plots of $\psi_{0,\alpha^\prime}^{\scriptstyle 2}(t)$   
and in Figure~\ref{fig.6} the
plots of $\psi_{1,\alpha^\prime}^{\scriptstyle 2}(t)$ for $\alpha^\prime = \{1,0.95,0.9\}$.

\section{Comparison between solutions of FDEs}

In this section we will compare the solutions of the equations with the 
dimensional regularized Caputo derivative obtained in the previous section with 
the solutions of the equation involving the regular Caputo derivative. 
For convenience we will choose $\omega = 1$.

\paragraph{(i) $\boldsymbol{N=1}$.} Let us consider the ODE 
\begin{equation}
\label{decay.order.1}
\frac{dy}{dt} +  y = 0 . 
\end{equation}
The corresponding FDE obtained using the regular Caputo derivative is 
\begin{equation}
\label{fde.ex.caputo.1}
{\sideset{_{}^{}}{_{{t}}^{\alpha}}{\operatorname{\mathcal{D}}}} y + y = 0 , \qquad 0 < \alpha \leq 1 ,
\end{equation}
whose solution is 
\begin{equation}
\label{sol.fde.ex.caputo.1}
y = y(0) \operatorname{E}_\alpha\left(- t^\alpha\right) ,
\end{equation}
where $\operatorname{E}_\alpha(\cdot)$ is the Mittag-Leffler function. The FDE obtained using
rule I or II is 
\begin{equation}
\label{fde.rule.1.bis}
t^{\alpha - 1} {\sideset{_{}^{}}{_{{t}}^{\alpha}}{\operatorname{\mathcal{D}}}}  y + y = 0 ,  \qquad 0 < \alpha \leq 1 ,
\end{equation}
whose solution is given by $\phi_{0,\alpha}^{\scriptstyle 1}( t)$ as in eq.\eqref{phi.zero.alpha}. 

\begin{figure}[hbt]
\begin{center}
\includegraphics[width=12cm]{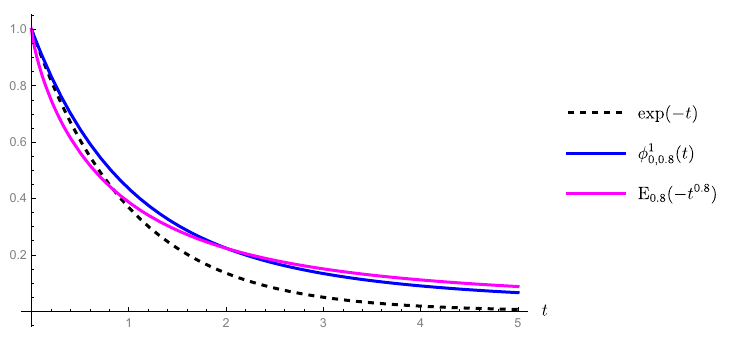}
\end{center}
\caption{Plots of the solution of eq.\eqref{decay.order.1}, that is, $\exp{(-t)}$, and its fractional version in eq.\eqref{fde.ex.caputo.1} and
eq.\eqref{fde.rule.1.bis} (with $\kappa = 1$)  for $\alpha = 0.8$, that is, $\operatorname{E}_{0.8}{(-t^{0.8})}$ and 
$\phi_{0,0.8}^{\scriptscriptstyle 1}(t)$, respectively. \label{fig.7}}
\end{figure}

In Figure~\ref{fig.7} we compare the solutions with $y(0) = 1$ of eq.\eqref{decay.order.1} and its fractional version in eq.\eqref{fde.ex.caputo.1} and
eq.\eqref{fde.rule.1.bis}. We choose for comparison $\alpha = 0.8$ and for simplicity we use $\kappa = 1$. So we   
have in Figure~\ref{fig.7} the plots of $\exp{(-t)}$, $\operatorname{E}_{0.8}{(-t^{0.8})}$ and 
$\phi_{0,0.8}^{\scriptstyle 1}(t)$. 

We can see that for values of $t$ such that $0 < t < 0.84$ (value obtained numerically using Mathematica), 
function $\operatorname{E}_{0.8}{(-t^{0.8})}$ decays faster than $\exp(-t)$, but for $t > 0.84$ this behaviour changes, 
and for $t \to \infty$ the function has a behaviour \cite{Mainardi}
$\operatorname{E}_{0.8}(-t^{0.8}) \sim t^{-0.8}/\Gamma(0.2)$.
The function $\phi_{0,0.8}^{\scriptscriptstyle 1}(t)$, on the other hand, has a decay that is always slower than $\exp(-t)$ for $t > 0$.
At the point $t = 2.01$ (value obtained numerically using Mathematica), the solution $\operatorname{E}_{0.8}(-t^{0.8})$ begins to 
show values greater than those of $\phi_{0,0.8}^{\scriptscriptstyle 1}(t)$, and after that these graphs never cross again. 
For $t\to \infty$ the we have $\phi_{0,0.8}^{\scriptscriptstyle 1}(t)  \sim t^{-1}/\Gamma(0.2)$ \cite{Simon}.

\paragraph{(ii) $\boldsymbol{N=2}$.} Let us consider the ODE  
\begin{equation}
\label{ode.two}
\frac{d^2 y}{dt^2} +  t = 0. 
\end{equation}
The corresponding FDE version obtained using the regular Caputo derivative is 
\begin{equation}
\label{fde.ex.caputo.2}
{\sideset{_{}^{}}{_{{t}}^{\alpha}}{\operatorname{\mathcal{D}}}} y +  y = 0 , \qquad 1 < \alpha \leq 2 ,
\end{equation}
whose solution is 
\begin{equation}
\label{sol.fde.ex.caputo.2}
y = y(0) \operatorname{E}_{\alpha,1}\left(- t^\alpha\right) + 
y^\prime(0) \, t \operatorname{E}_{\alpha,2}\left(- t^\alpha\right)  ,
\end{equation}
where $\operatorname{E}_{\alpha,\beta}(\cdot)$ is the two-parameter Mittag-Leffler function, 
\begin{equation}
\operatorname{E}_{\alpha,\beta}(z) = \sum_{n=0}^\infty \frac{z^n}{\Gamma(n\alpha + \beta)} ,
\end{equation}
and such that $\operatorname{E}_\alpha(\cdot) = \operatorname{E}_{\alpha,1}(\cdot)$ .
The FDE obtained using
rule I is 
 \begin{equation}
\label{fde.rule.10}
t^{\alpha - 2} {\sideset{_{}^{}}{_{{t}}^{\alpha}}{\operatorname{\mathcal{D}}}} y +  y  = 0 , 
\qquad 1 < \alpha \leq 2 , 
\end{equation}
whose solution is a linear combination $\phi^{\scriptstyle 2}_{0,\alpha}( t)$ 
and $\phi^{\scriptstyle 2}_{1,\alpha}(t)$ given by eq.\eqref{phi.zero.alpha.2} and
eq.\eqref{phi.one.alpha.2}, respectively. 
 On the other hand, the FDE obtained using rule II is
  \begin{equation}
\label{fde.rule.11}
\left(t^{\alpha^\prime - 1} {\sideset{_{}^{}}{_{{t}}^{\alpha^\prime}}{\operatorname{\mathcal{D}}}}\right)^2 y +  y  = 0 , 
\qquad 0 < \alpha^\prime \leq 1 , 
\end{equation}
whose solution is a linear combination $\psi^{\scriptstyle 2}_{0,\alpha^\prime}( t)$ 
and $\psi^{\scriptstyle 2}_{1,\alpha^\prime}(t)$ given by eq.\eqref{psi.0.alphaprime} and
eq.\eqref{psi.1.alphaprime}, respectively. 

\begin{figure}[hbt]
\begin{center}
\includegraphics[width=12cm]{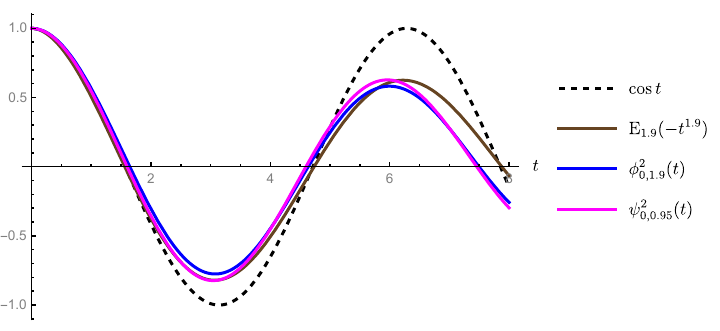}
\end{center}
\caption{Plots of the solutions of eq.\eqref{ode.two}, eq.\eqref{fde.ex.caputo.2}, eq.\eqref{fde.rule.10} and eq.\eqref{fde.rule.11} with initial conditions $y(0) =1$ and $y^\prime(0) = 0$.\label{fig.8}}
\end{figure}

Let us consider in separate the solutions for the initial conditions 
(i) $y(0) = 1$ and $ y^\prime(0) = 0$, and (ii) $y(0) = 0$ and $ y^\prime(0) = 1$. 
For (i) $y(0) = 1$ and $ y^\prime(0) = 0$, the solutions of 
eq.\eqref{ode.two}, eq.\eqref{fde.ex.caputo.2}, eq.\eqref{fde.rule.10} and 
eq.\eqref{fde.rule.11} are, respectively, $\cos{t}$, $\operatorname{E}_{1.9}\left(-t^{1.9}\right)$, 
$\phi_{0,1.9}^2(t)$ and $\psi_{0,0.95}^2(t)$, whose plots are in Figure~\ref{fig.8}.

\begin{figure}[hbt]
\begin{center}
\includegraphics[width=12cm]{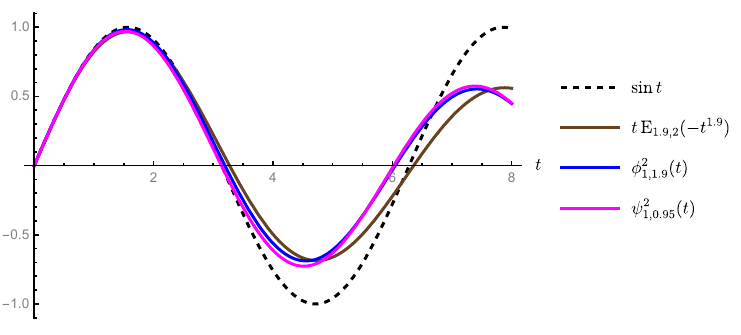}
\end{center}
\caption{Plots of the solutions of eq.\eqref{ode.two}, eq.\eqref{fde.ex.caputo.2}, eq.\eqref{fde.rule.10} and eq.\eqref{fde.rule.11} with initial conditions $y(0) =0$ and $y^\prime(0) = 1$.\label{fig.9}}
\end{figure}

For the initial conditions (ii) $y(0) = 0$ and $ y^\prime(0) = 1$, the solutions of 
eq.\eqref{ode.two}, eq.\eqref{fde.ex.caputo.2}, eq.\eqref{fde.rule.10} and 
eq.\eqref{fde.rule.11} are, respectively, $\sin{t}$, $t \operatorname{E}_{1.9}\left(-t^{1.9}\right)$, 
$\phi_{0,1.9}^2(t)$ and $\psi_{0,0.95}^2(t)$, whose plots are in Figure~\ref{fig.9}. 

Comparing the different plots in Figures~\ref{fig.8} and~\ref{fig.9}, 
we see that, as $t$ increases, the solutions of a FDE show a kind of damping, 
causing the amplitude of the functions $\operatorname{E}_{1.9}\left(-t^{1.9}\right)$, 
$\phi_{0,1.9}^2(t)$ and $\psi_{0,0.95}^2(t)$ in Figure~\ref{fig.8} and  of $t \operatorname{E}_{1.9}\left(-t^{1.9}\right)$, 
$\phi_{0,1.9}^2(t)$ and $\psi_{0,0.95}^2(t)$ in Figure~\ref{fig.9} to decrease. 
However, we note that these functions have a behaviour different than that of the solutions 
of the damped harmonic oscillator for large values of $t$. For the solutions of the damped harmonic oscillator, 
in addition to the decrease 
in the amplitude of the oscillations, these solutions also vanish at regularly spaced 
instants according to the period $\pi/\Omega$, where $\Omega$ is the damped frequency, 
and the oscillations continue until vanishing for $t\to \infty$. But this is not exactly 
what happens with the solutions of  eq.\eqref{fde.ex.caputo.2}, eq.\eqref{fde.rule.10} and eq.\eqref{fde.rule.11}. Indeed the intervals between successive instants where these solutions vanish 
are not regularly spaced like those of the solutions to the damped harmonic oscillator equation. 
In fact, for some large value $t_\ast$, although these solutions may still show a small oscillation, 
they will not vanish for $t > t_\ast$ because they tend to fit the asymptotic behaviour $\sim 1/t^\mu$ (for 
the appropriate value of $\mu$ in each case). 

\bigskip



\section{Concluding remarks}

We have seen that the process of writing a FDE by replacing an integer order derivative in a given DE by 
a fractional order derivative has some issues. The naive replacement of an integer order derivative 
by a fractional order derivative in a DE involving variables with physical dimensions leads to 
an equation with dimensional inconsistency. One way to try to get around this problem is to 
introduce parameters into the equation that compensate for the dimensional imbalance, 
but this introduces an element of arbitrariness into the FDE that is absent in the integer-order DE. 
On the other hand, an approach that does not need to resort to such parameters is to explore the 
dimensionless version of a DE. In fact, if the variables in a DE are dimensionless, the problem 
of dimensional balancing simply does not exist. In this case the naive replacement of an integer order 
derivative by a fractional order derivative does not bring any dimensional problems. 
However, in this case, when we reintroduce the physical dimensions into the fractional version of the DE, 
we have problems interpreting the terms of the FDE. We also have to make some choices involving 
characteristic values of the physical variables that introduce an element of arbitrariness into 
the process. But it is precisely the presence of this characteristic value with the same physical 
dimension as the independent variable of the DE that suggests to us how to define a version of Caputo  
fractional derivative that does not present the problems discussed. 
The dimensional regularized version of the Caputo fractional derivative we 
have discussed has the same physical dimension as the integer-order derivative 
for which it tends when $\alpha \to N$, and so the terms in the DE and in its
fractional version have the same interpretation. We solved some FDEs obtained
using two possible rules for using this dimensional regularization and we have
seen that the solutions can be conveniently written using the double gamma function 
and the $G$-Pochhammer symbol. We also compared these solutions with 
the ones for the analogous equation with the usual Caputo fractional derivative. 
The differences between these solutions were illustrated by some plots. 
Nevertheless these solutions present the same deviation characteristics from 
the solution of the integer order equation, such as a type of damping for the 
solutions of the equations of a harmonic oscillator, and also similar 
asymptotic behaviours.

\appendix

\section{The Barnes $\boldsymbol{G}$ and double gamma functions}
\label{appendix.A}

Around the turn of the 20th century, E. W. Barnes and others 
studied how to generalize the properties of the gamma function \cite{Barnes1,Barnes2,Barnes3,others}. 
Examples of functions defined in this context are the Barnes $G$-function $G(z)$ and 
its generalization $G(z,\tau)$, the so-called double gamma function, with 
$G(z,1) = G(z)$. Firstly let us 
define the function $G(z)$. 
The \textit{Barnes $G$-function} $G(z)$ was introduced in \cite{Barnes1} as a generalization 
of the gamma function $\Gamma(z)$. To see this, we begin by recalling the Weierstrass definition of the 
gamma function:  
\begin{equation}
\label{weierstrass.gamma}
\frac{1}{\Gamma(z)}  = z e^{\gamma z} \prod_{k=1}^\infty \left(1 + \frac{z}{k}\right) e^{-\frac{z}{k}} ,
\end{equation}
where $\gamma$ is the Euler-Mascheroni constant. This representation 
clearly shows that $1/\Gamma(z)$ is an entire function in the complex plane 
and that has simple zeros at $z = -k$ ($k=0,1,2,\ldots)$. 
Barnes defined the $G$-function as 
\begin{equation}
\label{weierstrass.G}
G(1+z)=(2\pi)^{z/2} \exp\left(- \frac{z+z^2(1+\gamma)}{2} \right) \, \prod_{k=1}^\infty \left\{ \left(1+\frac{z}{k}\right)^k \exp\left(\frac{z^2}{2k}-z\right) \right\} .
\end{equation}
Similar to $1/\Gamma(z)$, the function $G(z)$ is an entire function in the complex plane. 
However, the zeros of $G(z)$ at $z = -k$ have order  
$k+1$ at $z = -k$ ($k = 0,1,2,\ldots$). 

Equivalent definitions are 
\begin{equation}
\begin{split}
\label{ap.B.def.G.1}
G(z+1) = & (2\pi)^{z/2} \exp\left[\left(\gamma-1/2\right)z - \frac{z^2}{2}\left(
\frac{\pi^2}{6}+1+\gamma\right)\right] z \Gamma(z) \\
& \cdot \prod_{m,n\in \mathbb{N}_\ast^2} \left[ 
\left(1+\frac{z}{m+n}\right) \exp\left(-\frac{z}{m+n}+\frac{z^2}{2(m+n)^2}\right)\right] , 
\end{split}
\end{equation}
and 
\begin{equation}
\label{ap.B.def.G.2}
G(z+1)=(2\pi)^{z/2} \exp\left(- \frac{z+z^2(1+\gamma)}{2} \right) \, \prod_{k=1}^\infty \left\{
\frac{\Gamma(k)}{\Gamma(k+z)} \exp\left[z \psi(k) + \frac{z^2}{2}\psi^\prime(k)\right]\right\}
\end{equation}
where $\psi(\cdot)$ is the digamma function and $\mathbb{N}_\ast = \mathbb{N}\setminus\{0\}$. 

From these definitions we obtain the
functional relation 
\begin{equation}
\label{functional.barnes.G}
G(z+1) = \Gamma(z) G(z) . 
\end{equation}
Additionally, the definition implies that  
\begin{equation}
G(1) = 1. 
\end{equation}
As a result, we have 
\begin{equation}
G(n+2) = 0! \, 1! \, 2! \cdots n! \quad (n=0,1,2,\ldots) .
\end{equation}

 Denoted by $G(z; \tau)$, the \textit{double gamma function}  is a 
generalization of the $G$-function. 
It is defined as \cite{Kusnetsov}
\begin{equation}
\label{ap.B.G.tau.1}
\begin{split}
G(z;\tau) = & \frac{z}{\tau} \exp\left[a(\tau)\frac{z}{\tau} + b(\tau)\frac{z^2}{2\tau}\right] \\
& \cdot \prod_{m,n\in \mathbb{N}_\ast^2} \left[ 
\left(1+\frac{z}{m\tau+n}\right) \exp\left(-\frac{z}{m\tau+n}+\frac{z^2}{2(m\tau+n)^2}\right)\right] , 
\end{split}
\end{equation}
where $a(\tau)$ and $b(\tau)$ are given by 
\begin{equation}
\begin{split}
& a(\tau) = \gamma \tau + \frac{\tau}{2}\log{(2\pi \tau)} + \frac{1}{2}\log\tau - \tau C(\tau) , \\
& b(\tau) = - \frac{\pi^2 \tau^2}{6} - \tau \log\tau - \tau^2 D(\tau) , 
\end{split}
\end{equation}
with 
\begin{equation}
\begin{split}
& C(\tau) = \lim_{m\to \infty} \left[\sum_{k=1}^{m-1}\psi(k\tau) + \frac{1}{2}\psi(m\tau) - \frac{1}{\tau} 
\log\left(\frac{\Gamma(m\tau)}{\sqrt{2\pi}}\right) \right] , \\
& D(\tau) = \lim_{m\to \infty} \left[\sum_{k=1}^{m-1}\psi^\prime(k\tau) + \frac{1}{2}\psi^\prime(m\tau) - \frac{1}{\tau}\psi(m\tau) \right] . 
\end{split}
\end{equation}
Another definition of double gamma function is 
\begin{equation}
\label{ap.B.G.tau.2}
G(z;\tau) = \frac{1}{\tau\Gamma(z)} \exp\left[\tilde{a}(\tau)\frac{z}{\tau} + \tilde{b}(\tau) 
\frac{z^2}{2\tau^2}\right] \prod_{m=1}^\infty \frac{\Gamma(m\tau)}{\Gamma(z+m\tau)} 
\exp\left[z\psi(m\tau)+ \frac{z^2}{2}\psi^\prime(m\tau) \right] ,
\end{equation}
where 
\begin{equation}
\tilde{a}(\tau) = a(\tau) - \gamma\tau , \qquad \tilde{b}(\tau) = b(\tau) + \frac{\pi^2 \tau^2}{6} .
\end{equation}
Barnes $G$-function $G(z)$ is a particular case of the double 
gamma function $G(z;\tau)$ for $\tau = 1$. Note that $C(1) = 1/2$ and $D(1) = 1+ \gamma$, as shown 
in \cite{Barnes1}. 

The double gamma function is defined in such a way that 
\begin{equation}
\label{ap.B.G.tau.(1)}
G(1;\tau) = 1 . 
\end{equation}
It satisfies the functional relations 
\begin{equation}
\label{ap.B.general.G.tau}
G(z+1;\tau) = \Gamma(z/\tau) G(z;\tau) 
\end{equation} 
and 
\begin{equation}
\label{ap.B.general.G.tau.2}
G(z+\tau;\tau) = (2\pi)^{\frac{\tau-1}{2}} \tau^{-z+\frac{1}{2}} \Gamma(z) G(z;\tau) ,
\end{equation}
which generalize eq.\eqref{functional.barnes.G}. A straightforward consequence of these properties is 
\begin{equation}
\label{ap.B.G.tau.cons}
G(1+\tau;\tau) = G(\tau;\tau) = (2\pi)^{(\tau-1)/2} \tau^{-1/2} . 
\end{equation}
Moreover, using eq.\eqref{ap.B.general.G.tau}  and eq.\eqref{ap.B.general.G.tau.2}  recursively, we obtain
\begin{equation}
\label{ap.B.general.G.tau.n}
G(z+k;\tau) = G(z;\tau) \prod_{j=0}^{k-1} \Gamma[(z+j)/\tau]   
\end{equation}
and 
\begin{equation}
\label{ap.B.general.G.tau.2.n}
G(z+k\tau;\tau) = (2\pi)^{k(\tau-1)/2} \tau^{k(-z+(1-\tau(k-1))/2)} \left[\prod_{j=0}^{k-1} \Gamma(z+j\tau)\right] G(z;\tau)   , 
\end{equation}
respectively. 
Another interesting relations are  
\begin{equation}
\label{modular.transf}
G(z;\tau) = (2\pi)^{z(1-1/\tau)/2} 
\tau^{[(z-z^2)/(2\tau) +z/2 -1]} G(z/\tau;1/\tau) . 
\end{equation}
and, for $p, q \in \mathbb{N}$, 
\begin{equation}
\label{modular.2}
G(z;p\tau/q) = q^{\frac{1}{2p\tau}(z-1)(qz-p\tau)} (2\pi)^{-(q-1)(z-1)/2} 
\prod_{i=0}^{p-1} \prod_{j=0}^{q-1} \frac{G[(z+i)/p + j \tau/q;\tau]}{G(1+i)/p+j\tau/q;\tau]}  .
\end{equation}

An integral representation for $G(z;\tau)$ was provided in \cite{lawrie}, namely 
\begin{equation}
\label{double.gamma.int.rep}
\begin{aligned}
\ln G(z;\tau) & = \int_0^1 \bigg[ \frac{r^{z-1}}{(r-1)(r^\tau-1)} - 
\frac{z^2}{2\tau} r^{\tau-1} - z r^{\tau-1}\left(\frac{2-r^\tau}{r^\tau-1}-\frac{1}{2\tau}\right) \\
& -r^{\tau-1} + \frac{1}{r-1} - \frac{r^{\tau-1}}{(r-1)(r^\tau-1)}\bigg] \frac{dr}{\ln{r}} . 
\end{aligned}
\end{equation}
It is convergent for $\operatorname{Re}z > 0$ and $\delta > 0$.

Finally, let us recall that the Pochhammer symbol is defined as $(a)_n = a(a+1)\cdots (a+n-1)$ 
with $(a)_0 = 1$, or 
\begin{equation}
\label{ap.B.pochhammer}
(a)_n = \frac{\Gamma(a+n)}{\Gamma(a)} . 
\end{equation}
Since the double gamma function generalizes the gamma function, it is natural to generalize the 
definition of the Pochhammer symbol using $G(z;\tau)$ \cite{Simon}. So we define 
\begin{equation}
\label{ap.B.G.pochhammer}
\llbracket  a;\tau \rrbracket_n = \frac{G(a+n;\tau)}{G(a;\tau)} . 
\end{equation} 
From eq.\eqref{ap.B.general.G.tau.n} it follows that the usual Pochhammer symbol can be written in terms of the $G$-Pochhammer symbol 
as 
\begin{equation}
\label{ap.B.poch.G.poch.gen.2}
(a)_n = \tau^n \frac{\llbracket a+\tau;\tau\rrbracket_{n}}{\llbracket a;\tau\rrbracket_{n}} . 
\end{equation}

\section{The Kilbas-Saigo function}

Let $\alpha,m > 0$ and $l > -1/\alpha$. The Kilbas-Saigo function $E_{\alpha,m,l}(z)$ is defined by the series \cite{Meccanica} 
\begin{equation}
\label{KS.function.def}
E_{\alpha,m,l}(z) = \sum_{k=0}^\infty c_k z^k   
\end{equation}
with  
\begin{equation}
\label{def.ck}
c_0 = 1 , \qquad c_k = \prod_{j=0}^{k-1} \frac{\Gamma[\alpha(jm+l)+1]}{\Gamma[\alpha(jm+l+1)+1]} , \quad 
(k=1,2,\ldots) . 
\end{equation}

Let $\llbracket a;\tau\rrbracket_k$ be the $G$-Pochhammer symbol defined in eq.\eqref{ap.B.G.pochhammer}. 
Using eq.\eqref{ap.B.general.G.tau.n} it follows that the Kilbas-Saigo function $E_{\alpha,m,l}(z)$ can be written as 
\begin{equation}
\label{KS.G}
E_{\alpha,m,l}(z) = \sum_{k=0}^\infty \frac{\llbracket (\alpha l+1)/(\alpha m);1/(\alpha m)\rrbracket_k}{\llbracket (\alpha (l+1)+1)/(\alpha m);1/(\alpha m)\rrbracket_k} z^k .
\end{equation}

\end{document}